\documentclass[ejs,preprint]{imsart}

\RequirePackage[OT1]{fontenc}
\RequirePackage{amsthm,amsmath}
\RequirePackage[square,numbers,sort]{natbib}
\RequirePackage[colorlinks,citecolor=blue,urlcolor=blue]{hyperref}
\RequirePackage{comment}
\pubyear{0}
\volume{0}
\issue{0}
\firstpage{0}
\lastpage{0}
\startlocaldefs
\numberwithin{equation}{section}
\theoremstyle{plain}
\usepackage{amsmath,amssymb,amsfonts,amsthm,latexsym}
\usepackage{amscd,graphicx,color,enumerate,float,caption}
\usepackage[all]{xy}
\newtheorem{theorem}{Theorem}[section]
\newtheorem{lemma}[theorem]{Lemma}
\newtheorem{proposition}[theorem]{Proposition}

\theoremstyle{definition}
\newtheorem{definition}[theorem]{Definition}
\newtheorem{example}[theorem]{Example}
\theoremstyle{remark}
\newtheorem{remark}[theorem]{Remark}
\newcommand{\R}{\mathbb{R}}
\newcommand{\N}{\mathbb{N}}
\newcommand{\Prob}{\mathbb{P}}
\endlocaldefs
\begin{document}
%%%%%%%%%%%%%%%%%%%%%%%%%%%%%%%%%%%%%%%%%%%%%%%%%%%
\begin{frontmatter}
\title{The Kato--Temple inequality and eigenvalue concentration with applications to graph inference\thanksref{T1}}
\runtitle{Kato-Temple inequality, eigenvalue concentration, and graph inference}
\thankstext{T1}{This work is partially supported by the XDATA program of the Defense Advanced
	Research Projects Agency (DARPA) administered
	through Air Force Research Laboratory contract
	FA8750-12-2-0303, by the DARPA SIMPLEX
	program through SPAWAR contract N66001-15-C-4041,
	and by the DARPA GRAPHS program through contract N66001-14-1-4028. This work is also supported by the Acheson J. Duncan Fund for the Advancement of Research in Statistics at Johns Hopkins University.}

\begin{aug}
\author{\fnms{Joshua}
	\snm{Cape}\corref{}\ead[label=e1]{joshua.cape@jhu.edu}},
\author{\fnms{Minh}
	\snm{Tang}\ead[label=e2]{mtang10@jhu.edu}}, and
\author{\fnms{Carey E.}
	\snm{Priebe}\ead[label=e3]{cep@jhu.edu}}

\address{Department of Applied Mathematics and Statistics\\
	Johns Hopkins University\\
	3400 N. Charles St., Baltimore, MD 21218, USA
\printead{e1,e2,e3}}

\runauthor{J. Cape et al.}

\end{aug}

\begin{abstract}
We present an adaptation of the Kato--Temple inequality for bounding perturbations of eigenvalues with applications to statistical inference for random graphs, specifically hypothesis testing and change-point detection. We obtain explicit high-probability bounds for the individual distances between certain signal eigenvalues of a graph's adjacency matrix and the corresponding eigenvalues of the model's edge probability matrix, even when the latter eigenvalues have multiplicity. Our results extend more broadly to the perturbation of singular values in the presence of quite general random matrix noise.
\end{abstract}


%\tableofcontents
\end{frontmatter}
%%%%%%%%%%%%%%%%%%%%%%%%%%%%%%%%%%%%%%%%%%%%%%%%%%%
%%%%%%%%%%%%%%%%%%%%%%%%%%%%%%%%%%%%%%%%%%%%%%%%%%%
%%%%%%%%%%%%%%%%%%%%%%%%%%%%%%%%%%%%%%%%%%%%%%%%%%%
\section{Introduction}
\label{sec:Intro}
\subsection{Overview}
\label{sec:Overview}
Eigenvalues and eigenvectors are structurally fundamental quantities associated with matrices and are widely studied throughout mathematics, statistics, and engineering disciplines. For example, given an observed graph, the eigenvalues and eigenvectors of associated matrix representations (such as the adjacency matrix or Laplacian matrix) encode structural information about the graph (e.g. community structure, connectivity \cite{Chung1997}). In the context of certain random graph models, the eigenvalues and eigenvectors associated with the underlying matrix-valued model parameter, the edge probability matrix, exhibit similar information. It is therefore natural to study how ``close'' the eigenvalues and eigenvectors of a graph are to the underlying model quantities.

In this paper we consider simple, undirected random graphs on $n$ vertices generated via the \emph{inhomogeneous Erd\H{o}s--R\'{e}nyi model} (IERM) \cite{Hoff-et-al--2002,Bollobas2007}, $\mathbb{G}(n,P)$, where $P:=[P_{i,j}] \in [0,1]^{n \times n}$ denotes the (symmetric) edge probability matrix. This independent edge model generalizes numerous widely-studied random graph models including the classical Erd\H{o}s--R\'{e}nyi model \cite{Erdos-Renyi1959}, the stochastic block model \cite{Holland-et-al--1983}, and the random dot product graph model \cite{Young-Scheinerman--2007}. For $G \sim \mathbb{G}(n,P)$, the (symmetric) adjacency matrix, $A \equiv A_{G}\in\{0,1\}^{n \times n}$, has entries that are independently distributed according to $A_{i,j} \sim \text{Bernoulli}(P_{i,j})$ for all $i \le j$. This yields $P\equiv\mathbb{E}[A]$, where $\mathbb{E}[\cdot]$ denotes probabilistic expectation.

We focus our attention on the eigenvalues of $A$ and $P$. Specifically, we consider the eigenvalues in pairs (e.g. the largest eigenvalues of $A$ and of $P$ form a pair, as do the second-largest eigenvalues of each matrix, etc.). We obtain bounds on the distance between eigenvalues in certain ``signal pairs'', thereby and therein demonstrating a local sense in which random graphs concentrate. Note that in the random graph literature, the term \emph{concentration} is primarily used to describe global, uniform behavior via the spectral norm quantity $\|A-\mathbb{E}[A]\|_2$.

The following description provides an overview of our results for the IERM setting. Given a collection of consecutive, ordered eigenvalues of $P$ which are sufficiently separated from the remainder of the spectrum and conditional on that the corresponding eigenvalues of $A$ are not near the remainder of the spectrum of either $A$ or $P$, then Theorems \ref{thrm:IERMconditional} and \ref{thrm:spikeEvalsIERM} yield high-probability bounds on the distances between the eigenvalues in each pair. The individual, pair-specific (i.e. local) bounds we obtain stand in contrast to weaker bounds which hold uniformly for all eigenvalue pairs (for example, bounds implied by Weyl's inequality \cite{Horn-Johnson--2012}). Our results hold even in the presence of eigenvalue multiplicity.

We demonstrate that when the matrix $P$ has low rank, our results compare favorably with the recent study of low rank matrices undergoing random perturbation in \cite{O-Vu-Wang--2014} (see our Example \ref{ex:2BSSBM}). We also demonstrate that our results can lead to meaningful estimation in high rank settings (see Example \ref{ex:SpikeModel}).

After presenting our main theoretical results, we then apply the theory in this paper to both hypothesis testing and change-point detection for random graphs. Moreover, we generalize our results beyond the IERM setting to obtain high-probability bounds for perturbations of singular values of rectangular matrices in a quite general random matrix noise setting.

Broadly speaking, we adapt the original, deterministic setting in a paper by T. Kato \cite{Kato--1950} to a new setting involving randomness, and this approach is novel in the context of random graphs, random matrix theory, and statistical inference for random graphs. We further detail the key modifications and differences between our work and \cite{Kato--1950} in our subsequent remarks and proofs. The present paper also stands in contrast to a deterministic generalization of the Kato--Temple inequality in \cite{Harrell--1978}.

\subsection{Inhomogeneous random graphs}
In the inhomogeneous random graph literature, concentration bounds have been known for some time for each eigenvalue of $A$, denoted $\lambda_{i}(A)$, both around its median and around its expectation, $\mathbb{E}[\lambda_{i}(A)]$ \cite{Alon-K-Vu--2002}. Unfortunately, since the latter quantities are inaccessible in practice, such bounds are of limited practical use. Moreover, in general $\mathbb{E}[\lambda_{i}(A)] \neq \lambda_{i}(\mathbb{E}[A])$.

By way of contrast, numerous results in the literature bound the spectral norm matrix difference $\|A-\mathbb{E}[A]\|_{2}$, thereby immediately and uniformly bounding each of the eigenvalue differences $|\lambda_{i}(A)-\lambda_{i}(\mathbb{E}[A])|$ via an application of Weyl's inequality. For example, \cite{Oliveira--2010} proved an asymptotically almost surely spectral norm bound of $\|A-\mathbb{E}[A]\|_2 =  O(\sqrt{\Delta \log n})$ for $\Delta=\Omega(\log(n))$ where $\Delta\equiv\Delta(n)$ denotes the maximum expected degree of a graph.
In \cite{Lu-Peng--2013} the above bound is improved to $\|A-\mathbb{E}[A]\|_2 \le (2+o(1))\sqrt{\Delta}$ under the stronger assumption that $\Delta = \omega(\log^4 n)$ with further refinement being subsequently obtained in \cite{LeiRinaldo2015}. We on the other hand show that under certain conditions, for particular eigenvalue pairs one can obtain tighter and non-uniform high probability bounds of the form $|\lambda_{i}(A)-\lambda_{i}(\mathbb{E}[A])|=O(\log^{\delta}n)$ for small $\delta>0$.

Spectral theory for random graphs overlaps with the random matrix theory literature.
%insofar as graphs may be studied as matrices
There, asymptotic analysis includes proving, for example, convergence of the empirical spectral distribution to a limiting measure \cite{Ding-Tiefeng2010}. Related approaches to studying the spectrum of random graphs consider normalized versions of the adjacency matrix \cite{Le-Vershynin2015} and employ standard random matrix theory techniques such as the Stieltjes transform method \cite{Avrachenkov-Cottatellucci-Kadavankandy2015, Zhang-Nad.-Rao--2014}. In contrast, we do not study normalized versions of the adjacency or the edge probability matrix.

Indeed, much of the existing literature focuses on properties of eigenvectors corresponding to random graphs \cite{Fortunato2010,LeiRinaldo2015,TangPriebe2016} given, among other reasons, the success of spectral clustering methods for graph inference \cite{Luxburg2007}. We do not consider eigenvectors since our aim is to demonstrate the usefulness of adapting and applying the eigenvalue-centric Kato--Temple framework.

The stochastic block model (SBM) offers an example of an inhomogeneous random graph model which is wildly popular in the literature \cite{Lei2016,BickelSarcar2013,Zhao2012,KarrerNewman2011,LeiRinaldo2015} and in which our results apply to the top (signal) eigenvalues of $A$ and $P$. Previously, the authors in \cite{Athreya-et-al--2015} obtained a collective deviation bound on the top eigenvalues of $A$ and $P$ for certain stochastic block model graphs in order to prove the main limit theorem therein. Our Theorem \ref{thrm:IERMconditional} improves upon Lemma 2 in \cite{Athreya-et-al--2015} by removing a distinct eigenvalue assumption and by yielding stronger high-probability deviation bounds for pairs of top eigenvalues of $A$ and $P$ which are of the same order. This implies a statistical hypothesis testing regime for random graphs which is discussed further in Section \ref{sec:Application}.
%%%%%%%%%%%%%%%%%%%%%%%%%%%%%%%%%%%%%%%%%%%%%%%%%%%
\subsection{Organization}
\label{sec:Organization}
The remainder of this paper is organized as follows.
In Section \ref{sec:Setup} we introduce notation and the Kato--Temple eigenvalue perturbation framework.
In Section \ref{sec:Results} we present our results for random graphs and more generally for matrix perturbation theory. There we also include illustrative examples together with comparative analysis involving recent results in the literature. 
In Section \ref{sec:Application} we discuss applications of our results to problems involving graph inference.
Sections \ref{sec:Thanks} and \ref{sec:Appendix} contain our acknowledgments and the proofs of our results, respectively.
%%%%%%%%%%%%%%%%%%%%%%%%%%%%%%%%%%%%%%%%%%%%%%%%%%%
%%%%%%%%%%%%%%%%%%%%%%%%%%%%%%%%%%%%%%%%%%%%%%%%%%%
%%%%%%%%%%%%%%%%%%%%%%%%%%%%%%%%%%%%%%%%%%%%%%%%%%%
\section{Setup and notation}
\label{sec:Setup}
Let $\langle \cdot, \cdot \rangle$ denote the standard Euclidean inner (dot) product between two vectors, $\|\cdot\|$ denote the vector norm induced by the dot product, and $\|\cdot\|_{2}$ denote the spectral norm of a matrix. The identity matrix is implicitly understood when we write the difference of a matrix with a scalar. In this paper, $\mathcal{O}(\cdot)$, $\Omega(\cdot)$, and $\Theta(\cdot)$ denote standard big-O, big-Omega, and big-Theta notation, respectively, while $o(\cdot)$ and $\omega(\cdot)$ denote standard little-o and little-omega notation, respectively.

As prefaced in Section \ref{sec:Intro}, we consider simple, undirected random graphs on $n$ vertices generated by the inhomogeneous Erd\H{o}s--R\'{e}nyi model, $G \sim \mathbb{G}(n,P)$, via the corresponding (binary, symmetric) adjacency matrix $A \equiv A_G$.
Given an open interval in the positive half of the real line, $(\alpha, \beta) \subset \R_{>0}$, we denote the $d$ eigenvalues of $P$ that lie in this interval (locally) by
\begin{equation}
	\label{eqn:local}
	\alpha
	< \lambda_{1}(P)
	\le \lambda_{2}(P)
	\le \dots
	\le \lambda_{d}(P)
	< \beta,
\end{equation}
and similarly for $A$, noting that for $A$ this amounts to a probabilistic statement. By symmetry one can just as well handle the case when the interval lies in the negative half of the real line. We are principally interested in eigenvalues that are large in magnitude, so we do not consider the case when the underlying interval contains the origin.

To highlight the Kato--Temple framework for bounding eigenvalues, we now reproduce two lemmas from \cite{Kato--1950} along with the Kato--Temple inequality as stated in \cite{Harrell--1978} (see Theorem \ref{thrm:KT_single_case} below).\footnote{Of primary importance in this paper is the extension of Theorem \ref{thrm:KT_single_case} to multiple eigenvalues as presented in \cite{Kato--1950}. The original statement of the extension to multiple eigenvalues is more involved and therefore omitted for simplicity.} These results all hold in the following common setting.

\begin{quote}
	Let $H$ be a self-adjoint operator on a Hilbert space. Assume a unit vector $w$ is in the domain of $H$ and define $\eta:=\langle Hw,w \rangle$ along with $\epsilon := \|(H-\eta)w\|$, noting that $\eta^2 + \epsilon^2 = \|Hw\|^2$. The quantity $\eta$ may be viewed as an ``approximate eigenvalue'' of $H$ corresponding to the ``approximate eigenvector'' $w$, while $\epsilon$ represents a scalar residual term.
\end{quote}

\begin{lemma}[\cite{Kato--1950}, Lemma 1]
	\label{lem: kato1}
	For every $\alpha$ such that $\alpha < \eta$ (where $\alpha = -\infty$ is permitted), the interval $(\alpha, \eta + \frac{\epsilon^2}{\eta - \alpha}]$ contains a point in the spectrum of $H$.
\end{lemma}

\begin{lemma}[\cite{Kato--1950}, Lemma 2]
	\label{lem:kato2}
	For every $\beta$ such that $\beta > \eta$ (where $\beta = \infty$ is permitted), the interval $[\eta - \frac{\epsilon^2}{\beta - \eta}, \beta)$ contains a point in the spectrum of $H$.
\end{lemma}

\begin{theorem}[Kato--Temple inequality; \cite{Harrell--1978}, Theorem 2]
	\label{thrm:KT_single_case}
	Suppose that $\epsilon^2 < (\beta - \eta)(\eta - \alpha)$ where $\alpha < \beta$. Then $\emph{spectrum}(H) \cap (\alpha, \beta) \neq \emptyset$. Moreover, if the only point of the spectrum of $H$ in the interval $(\alpha, \beta)$ is the eigenvalue $\lambda(H)$, then
	$$ - \frac{\epsilon^2}{\beta - \eta}
	\le \lambda(H) - \eta
	\le \frac{\epsilon^2}{\eta - \alpha}.$$
\end{theorem}

\begin{remark}[Hermitian dilation]
	\label{rem:HermDial}
	Given an $m \times n$ real matrix $M$, it will be useful to consider the corresponding real symmetric $(m+n) \times (m+n)$ Hermitian dilation matrix $\tilde{M}$ given by
	\[
	\tilde{M} := \left[ \begin{array}{cc}
	0 & M \\
	M^\top & 0 \end{array} \right].
	\]
	It is well-known that the non-zero eigenvalues of $\tilde{M}$ correspond to the signed singular values of $M$ (see Theorem 7.3.3 in \cite{Horn-Johnson--2012}).
	This correspondence between the singular values of arbitrary matrices and the eigenvalues of Hermitian matrices allows our results to generalize beyond the IERM setting to the more general study of matrix perturbation theory for singular values in a straightforward manner.
\end{remark}
%%%%%%%%%%%%%%%%%%%%%%%%%%%%%%%%%%%%%%%%%%%%%%%%%%%
%%%%%%%%%%%%%%%%%%%%%%%%%%%%%%%%%%%%%%%%%%%%%%%%%%%
%%%%%%%%%%%%%%%%%%%%%%%%%%%%%%%%%%%%%%%%%%%%%%%%%%%
\section{Results}
\label{sec:Results}
%%%%%%%%%%%%%%%%%%%%%%%%%%%%%%%%%%%%%%%%%%%%%%%%%%%
\subsection{Results for random graphs}
\label{sec:ResultsRandomGraphs}
In the IERM setting, a graph's adjacency matrix can be written as $A=P+E$ where $E:=A-P$ is a random matrix and $P$ is the (deterministic) expectation of $A$. We begin with a preliminary observation concerning the tail behavior of $A-P$ which will subsequently be invoked for the purpose of obtaining standard union bounds. The proof follows from a straightforward application of Hoeffding's inequality.
\begin{proposition}[General IERM concentration]
	\label{prop:A-P}
	Let $u, v \in\R^{n}$ denote (non-random) unit vectors. Then for any $t>0$,
	\begin{equation}
		\label{eqn:RDPG tail bound}
		\Prob[|\langle(A-P)u,v\rangle| > t]
		\le 2 \exp(-t^2).
	\end{equation}
\end{proposition}
It is indeed possible to invoke more refined concentration inequalities than Proposition \ref{prop:A-P} in the presence of additional structure (e.g. when all entries of $P$ have uniformly very small magnitude). Doing so is particularly useful when it is simultaneously possible to obtain a strong bound on $\|A-P\|_{2}$. This observation will be made clearer in the context of Theorem \ref{thrm:IERMconditional} below. Furthermore, consideration of Proposition \ref{prop:A-P} will facilitate the subsequent presentation of our generalized results which extend beyond the IERM setting.
\begin{remark}
	\label{rem:diagElementsOfP}
	In this paper the main diagonal elements of $P$ are allowed to be strictly positive, in which case realizations of $A$ need not necessarily be hollow (i.e. observed graphs may have self-loops). To avoid graphs with self-loops, one may either condition on the event that $A$ is hollow or set the main diagonal of $P$ to be zero. In the former case, note that $P\equiv\mathbb{E}[A]$ no longer holds on the main diagonal. In the latter case, a modified version of Proposition \ref{prop:A-P} holds.
\end{remark}

\begin{comment}
\begin{remark}
	\label{rem:A-Ptilde}
	For each $i \in [n]$ define the vector $\tilde{w}_i := [\frac{1}{\sqrt{2}}w_i, \frac{1}{\sqrt{2}}w_i]^\top \in \R^{2n}$. Observe that by construction $\langle \tilde{A}\tilde{w}_i, \tilde{w}_j \rangle
	= \langle A w_i, w_j \rangle$ for all $i,j \in [n]$.
	%	since $$ \tilde{w}_j^T \tilde{A} \tilde{w}_i
	%	= (\frac{1}{\sqrt{2}}w_j)^T A (\frac{1}{\sqrt{2}}w_i) + (\frac{1}{\sqrt{2}}w_j)^T A^T (\frac{1}{\sqrt{2}}w_i)
	%	= w_j^T A w_i.$$
	Moreover, the same relationship holds for matrices $\tilde{P}$ and $\tilde{A}-\tilde{P}$. Hence, Proposition \ref{prop:A-P} immediately extends to the collection $\tilde{A}, \tilde{P}$, and $\{\tilde{w}_i\}_{i=1}^{n}$.
\end{remark}
\end{comment}

We now present our main results for the IERM setting. The proofs, which are located in Section \ref{sec:Appendix}, also formulate a bound for the special case when the upper bound threshold $\beta$ may be chosen to be infinity. This special case is particularly useful in applications.

\begin{theorem}[IERM eigenvalue perturbation bounds, conditional version]
	\label{thrm:IERMconditional}
	Let the matrices $A \in \{0,1\}^{n \times n}$ and $P \in [0,1]^{n \times n}$ correspond to the IERM setting described in Section $\ref{sec:Setup}$. Suppose the interval $(\alpha,\beta) \subset \R_{>0}$ contains precisely $d$ eigenvalues of $P$, $\lambda_{1}(P)\le\lambda_{2}(P)\le \dots \le \lambda_{d}(P)$ (possibly with multiplicity). Condition on the event that $(\alpha, \beta)$ contains precisely $d$ eigenvalues of $A$, $\{\lambda_{i}(A)\}_{i=1}^{d}$, as well as the set $\{\langle A w_{i}, w_{i} \rangle\}_{i=1}^{d}$ where $\{w_i\}_{i=1}^{d}$ is an orthonormal collection of eigenvectors of $P$ corresponding to the eigenvalues $\{\lambda_{i}(P)\}_{i=1}^{d}$. Fix $k \in [d]$. Define $l:=(d-k+1)$. Then, for $t>0$,
	\begin{align}
		\label{eqn:LowerBound}
		\lambda_{k}(A)
		\ge \lambda_{k}(P) - t
		&- \zeta^{-},
	\end{align}
	where
	$\zeta^{-} := \frac{l\|E\|_{2}^{2}+((\beta-\lambda_{k}(P))+(\lambda_{d}(P)-\lambda_{k}(P)) + 3t)l(l-1)t}{\beta-\lambda_{d}(P)-(l(l-1)+1)t}$
	with probability at least
	$1-\left(l+\binom{l}{2}\right)2\exp(-t^{2})$.
	Also, for $t>0$,
	\begin{align}
		\label{eqn:UpperBound}
		\lambda_{k}(A)
		\le \lambda_{k}(P) + t
		&+ \zeta^{+},
	\end{align}
	where $\zeta^{+}:=\frac{k\|E\|_{2}^{2}+(3\lambda_{k}(P)-\alpha +3t)k(k-1)t}{\lambda_{1}(P)-\alpha - (k(k-1)+1)t}$
	with probability at least\\
	$1-\left(k+\binom{k}{2}\right)2\exp(-t^{2})$.
	Moreover, the upper and lower bounds hold collectively with probability at least $1-\left(d+\binom{d}{2}\right)2\exp(-t^{2})$.
\end{theorem}

\begin{remark}
	Our proof depends upon several new observations with respect to Kato's original argument. In particular, for $w_i$ as defined above, the matrix $[\langle A w_i, w_j \rangle]_{i,j=1}^{d}$ need not be diagonal, so $\{w_i\}_{i=1}^{d}$ need not constitute an orthonormal collection of ``approximate eigenvectors'' of $A$ in the sense of \cite{Kato--1950}. Instead, here the notion of ``approximate'' may be interpreted via Proposition \ref{prop:A-P} as the source of randomness which allows for Kato--Temple methodology to be adapted beyond the original deterministic setting. Of additional note is that the vectors $w_i$ as defined in this paper agree in function and notation with Kato's original paper, the operational distinction being that our setting provides a canonical choice for these vectors.
\end{remark}
\begin{remark}
	\label{rem:spec_vs_Euclid bounds}
	We note that the term $\|E\|_{2}^{2}$ in the formulation of both $\zeta^{+}$ and $\zeta^{-}$ can be replaced by an appropriate maximum over quantities of the form $\|Ew_{i}\|^{2}$ (see Equation (\ref{eqn:epsilonSquared})). That is to say, in the presence of additional local structure and knowledge, one can refine the above bounds in Theorem \ref{thrm:IERMconditional}.
\end{remark}
\begin{remark}
	\label{rem:interval extension}
	In settings wherein the eigenvalues of interest have disparate orders of magnitude, Kato--Temple methodology is not guaranteed to yield useful bounds. This can be seen in the bounds' dependence on the ratio of the eigenvalues of $P$ in Theorem \ref{thrm:IERMconditional}. Moreover, within the Kato--Temple framework, poor separation from the remainder of the spectrum also deteriorates the bounds, as is evident in the denominators' dependence on the interval endpoints  $\alpha$ and $\beta$ along with the smallest and largest local eigenvalues of $P$. On the other hand, by further localizing, i.e. by restricting to a subset of $d^{\prime} < d$ eigenvalues in a particular interval, applying Theorem \ref{thrm:IERMconditional} to said fewer eigenvalue pairs may yield improved bounds (see Example \ref{ex:2BSSBM} and Remark \ref{rem:spec_vs_Euclid bounds}).
\end{remark}

Next, we formulate an unconditional version of Theorem \ref{thrm:IERMconditional}. For both simplicity and the purpose of applications, Theorem \ref{thrm:spikeEvalsIERM} is stated in terms of the largest singular values in the IERM setting.

\begin{theorem}[IERM singular value perturbation bounds, unconditional version]
	\label{thrm:spikeEvalsIERM}
	Let the matrices $A \in \{0,1\}^{n \times n}$ and $P \in [0,1]^{n \times n}$ correspond to the IERM setting described in Section $\ref{sec:Setup}$ with maximum expected degree (via $P$) given by $\Delta\equiv\Delta(n)$. Denote the $d+1$ largest singular values of $A$ by $0 \le \hat{\sigma}_{0} < \hat{\sigma}_{1} \le \dots \le \hat{\sigma}_{d}$, and denote the $d+1$ largest singular values of $P$ by $0 \le \sigma_{0} < \sigma_{1} \le \dots \le \sigma_{d}$. Suppose that
	$\Delta =\omega(\log^{4}n)$,
	$\sigma_{1} \ge C\Delta$,
	and  $\sigma_{0} \le c\Delta$
	for some absolute constants $C>c>0$.
	Let $\delta\in(0,1]$. Then for each $k\in[d]$, there exists some positive constant $c_{k,d}$ such that as $n \rightarrow\infty$, with probability $1-o(1)$ involving $\delta$,
	\begin{equation}
		|\hat{\sigma}_{k}-\sigma_{k}|\le c_{k,d}\left(\log^{\delta}n\right).
	\end{equation}
\end{theorem}
A similar version of Theorem \ref{thrm:generalSVbound} holds when $\Delta=\Omega(\log n)$ under slightly different assumptions on the entries of $P$ for which one still has $\|A-P\|_{2}=O(\sqrt{\Delta})$ with high probability \cite{LeiRinaldo2015}. On a related yet different note, see \cite{Le-Vershynin2015} for discussion of the sparsity regime $\Delta=O(1)$ in which graphs fail to concentrate in the classical sense.
\begin{remark}[Random dot product graph model]
	\label{rem:RDPG}
	When the edge probability matrix $P$ can be written as $P=XX^{\top}$ for some matrix $X\in\R^{n \times d}$ with $d \ll n$, then the IERM corresponds to the popular \emph{random dot product graph (RDPG) model} \cite{Young-Scheinerman--2007}. In the random dot product graph model, the largest eigenvalues of $A$ and $P$ are of statistical interest in that they represent spectral ``signal'' in the model. These eigenvalues are separated from the remainder of their respective spectra and lie in an interval of the form $(\alpha, \infty)$ where, for example, $\alpha$ may be taken to be $O(\|A-P\|_2)$.
	
	Among its applications, the RDPG model has been used as a platform for modeling graphs with hierarchical and community structure \cite{lyzinski2017community}. A central limit theorem is known for the behavior of the top eigenvectors of adjacency matrices arising from the RDPG model \cite{Athreya-et-al--2015}, and this limit theorem relies upon a lemma which collectively bounds the differences between top eigenvalues of $A$ and $P$ but requires a stringent eigengap assumption.
	Namely, Lemma 2 in \cite{Athreya-et-al--2015} states that with high probability,
	\begin{equation}
		\label{eqn:CLT big 0_P bound}
		\sqrt{\sum_{i=1}^{d}|\lambda_{i}(A) - \lambda_{i}(P)|^2} = O(\delta_{\textnormal{gap}}^{-2} \log{n}).
	\end{equation}
	In contrast, using Theorem \ref{thrm:spikeEvalsIERM} with $\sigma_{0}:=0$, we do not require the gap assumption $\delta_{\text{gap}} > 0$ and still obtain that with high probability,
	\begin{equation}
		\label{eqn:CLT_Kato_improved}
		\sqrt{\sum_{i=1}^{d}|\lambda_{i}(A) - \lambda_{i}(P)|^2} = O(\log{n}).
	\end{equation}
	In practice, models involving repeated or arbitrarily close eigenvalues are prevalent and of interest (e.g. Section \ref{subsec:threeBlockSBM}). As such, the above improvement is nontrivial and of practical significance.
\end{remark}

\begin{remark}[Latent position random graphs]
	\label{rem:kernel extension}
	Theorem \ref{thrm:IERMconditional} further extends to the more general setting of latent position random graphs. There, the matrix $P$ is viewed as an operator $[\kappa(X_i, X_j)]_{i,j=1}^{n}$ where $X_{i}$ and $X_{j}$ are independent, identically distributed latent positions with distribution $F$, and the positive definite kernel, $\kappa$ (viewed as an integral operator), is not necessarily of finite fixed rank as $n$ increases \cite{Hoff-et-al--2002, Tang-et-al--2013}. Note that for the RDPG model, the kernel $\kappa$ is simply the standard Euclidean inner product between (latent position) vectors.
\end{remark}

%%%%%%%%%%%%%%%%%%%%%%%%%%%%%%%%%%%%%%%%%%%%%%%%%%%
\subsection{Results for matrix perturbation theory}
\label{sec:ResultsExtension}
The behavior of the random matrix  $A-P$ (see Proposition \ref{prop:A-P}) represents a specific instance of more general, widely-encountered probabilistic concentration as discussed in \cite{O-Vu-Wang--2014} and formulated in the following definition.

\begin{definition}[\cite{O-Vu-Wang--2014}]
	\label{def:Ccgamma}
	An $m \times n$ random real (``error'') matrix $E$ is said to be \emph{$(C,c,\gamma)$-concentrated} for a trio of positive constants $C,c,\gamma >0$ if for all unit vectors $u \in \R^{n}, v \in \R^{m}$ and for every $t>0$, then
	\begin{equation}
	\label{eqn:concentration}
	\mathbb{P}[|\langle Eu, v \rangle| > t]
	\le C \exp(-c t^\gamma).
	\end{equation}
\end{definition}

In particular, the IERM setting corresponds to $(C,c,\gamma)$ concentration where $m=n$, $C=\gamma=2$, and $c=1$. For the Hermitian dilation discussed in Remark \ref{rem:HermDial} one has the following correspondence between $E$ and $\tilde{E}$.

\begin{lemma}[\cite{O-Vu-Wang--2014}]
	\label{lem:CcgammaBlowUp}
	Let $E\in\R^{m\times n}$ be $(C,c,\gamma)$-concentrated. Define $\tilde{C}:=2C$ and $\tilde{c}:=c/2^{\gamma}$. Then the matrix $\tilde{E}\in\R^{m+n \times m+n}$ is $(\tilde{C},\tilde{c},\gamma)$-concentrated.
\end{lemma}

Definition \ref{def:Ccgamma} and Lemma \ref{lem:CcgammaBlowUp} together with Remark \ref{rem:HermDial} allow for Theorem \ref{thrm:IERMconditional} to be generalized in a straightforward manner. We frame the generalization in the context of a signal--plus--noise matrix model with tail probability bounds. In particular, replace $A$ with $\hat{M}:=M+E$, thought of as an observed data matrix. Also replace $P$ with $M$, thought of as an underlying signal matrix, so that the matrix $A-P$ becomes $E$, thought of as an additive error matrix. We emphasize that the following generalization is in terms of the singular values of $M$ and $\hat{M}$. This generalization resembles the formulation of a result obtained in \cite{O-Vu-Wang--2014} using different methods; however, unlike our Theorem \ref{thrm:generalSVbound}, the bound in \cite{O-Vu-Wang--2014} depends upon the rank of $M$ and assumes that the rank is known.

Given a matrix $M\in\R^{m\times n}$, write its singular value decomposition as $M\equiv U\Sigma V^{\top}$ where $Mv_{i}=\sigma_{i}u_{i}$ holds for the normalized left (resp., right) singular vectors $v_{i}$ (resp., $u_{i}$) and singular values $\sigma_{i}=\Sigma_{i,i}$.
For each $i$ such that $\sigma_{i}>0$, define $\tilde{w}_{i}\in\R^{m+n}$ to be the concatenated unit vector $w_{i}:=\frac{1}{\sqrt{2}}(u_{i}^{\top},v_{i}^{\top})^{\top}$. Note that $\tilde{w}_{i}$ is an eigenvector for $\tilde{M}$ with $\tilde{M}\tilde{w}_{i}=\sigma_{i}\tilde{w}_{i}$.
\newpage
\begin{theorem}[Singular value perturbation bounds, conditional version]
	\label{thrm:generalSVbound}
	For matrices $M,E \in \R^{m \times n}$ and $\hat{M}:=M+E$, suppose that $E$ is ($C,c,\gamma$)-concentrated for positive constants $C, c, \gamma >0$. Suppose the interval $(\alpha,\beta)\subset\R_{>0}$ contains the largest $d$ singular values of $M$, denoted by 
	$0<\sigma_{1}\le\sigma_{2}\le\dots\le\sigma_{d}$. Condition on the event that the interval $(\alpha,\beta)$ contains precisely $d$ singular values of $\hat{M}$, denoted 
	$0<\hat{\sigma}_{1}\le\hat{\sigma}_{2}\le\dots\le\hat{\sigma}_{d}$, as well as
	$\langle \tilde{\hat{M}}\tilde{w}_{i},\tilde{w}_{i}\rangle$ for $1 \le i \le d$ and unit vector $\tilde{w}_{i}$ as defined above.
	Fix $k \in [d]$. Define $l:=(d-k+1)$.
	Then for $t > 0$,
	\begin{align}
		\label{eqn:GeneralLowerBound}
		\hat{\sigma}_{k}
		\ge \sigma_{k} - t
		&- \zeta^{-},
	\end{align}
	where
	$\zeta^{-} := \frac{l\|E\|_{2}^{2}+((\beta-\sigma_{k})+(\sigma_{d}-\sigma_{k}) + 3t)l(l-1)t}{\beta-\sigma_{d}-(l(l-1)+1)t}$
	with probability at least\\
	$1-\left(l+\binom{l}{2}\right)\tilde{C}\exp(-\tilde{c}t^{\gamma})$.
	Also, for $t>0$,
	\begin{align}
		\label{eqn:GeneralUpperBound}
		\hat{\sigma}_{k}
		\le \sigma_{k} + t
	&+ \zeta^{+},
	\end{align}
	where $\zeta^{+}:=\frac{k\|E\|_{2}^{2}+(3\sigma_{k}-\alpha +3t)k(k-1)t}{\sigma_{1}-\alpha - (k(k-1)+1)t}$
	with probability at least\\
	$1-\left(k+\binom{k}{2}\right)\tilde{C}\exp(-\tilde{c}t^{\gamma})$. Moreover, the upper and lower bound hold collectively with probability at least $1-\left(d+\binom{d}{2}\right)\tilde{C}\exp(-\tilde{c}t^{\gamma})$.
\end{theorem}

As with the results in Section \ref{sec:ResultsRandomGraphs}, Theorem \ref{thrm:generalSVbound} can be formulated unconditionally and for collections of not--necessarily--the--largest singular values. Both of these aspects are explored in greater detail in Example \ref{ex:SpikeModel}. The following technical lemma will subsequently be employed in applications.

\begin{lemma}
	\label{lem:CcgammaSpectralProb}
	Let $E\in\R^{m \times n}$ be a $(C,c,\gamma)$-concentrated random matrix.
	Choose $\epsilon > 0$ such that $2 + \epsilon > 2\left(2\log(9)/c\right)^{1/\gamma}$ and define the quantity $c_{\epsilon,c,\gamma}
	:=\left(c(1+\epsilon/2)^{\gamma}-2\log(9)\right)>0$. Then,
	\begin{align}
	\mathbb{P}\left[\|E\|_{2}>(2+\epsilon)\textnormal{max}\{m,n\}^{1/\gamma}\right]
	&\le C\exp(-c_{\epsilon,c,\gamma}\textnormal{max}\{m,n\}).
	\end{align}
	If in addition $m=n$ and $E$ is assumed to be symmetric, then the quantity $2\log(9)$ above may be replaced by $\log(9)$, an improvement.
\end{lemma}
%%%%%%%%%%%%%%%%%%%%%%%%%%%%%%%%%%%%%%%%%%%%%%%%%%%
%%%%%%%%%%%%%%%%%%%%%%%%%%%%%%%%%%%%%%%%%%%%%%%%%%%
%%%%%%%%%%%%%%%%%%%%%%%%%%%%%%%%%%%%%%%%%%%%%%%%%%%
\subsection{Two illustrative examples}
\label{sec:Discussion}
In the remainder of this section, we present two examples that highlight the usefulness and flexibility of Kato--Temple methodology. We begin with Example \ref{ex:2BSSBM} which presents a simple stochastic block model setting wherein our results compare favorably with those in the recent work of \cite{O-Vu-Wang--2014}, noting that in general for similar settings, the corresponding results are often comparable.
%%%%%%%%%%%%%%%
\begin{example}[Balanced affinity two block stochastic block model]
	\label{ex:2BSSBM}
	Consider an $n$ vertex realization from a two block stochastic block model in which $0<q<p<1$ where $p$ and $q$ denote the within--block and between--block edge probabilities, respectively. Suppose each block contains $n/2$ of the graph's vertices. The signal singular values and maximum expected degree of this rank two model are given by
	\begin{equation}
		\sigma_{1}(P)=\frac{n}{2}(p-q),
		\sigma_{2}(P)=\frac{n}{2}(p+q), \textnormal{ and }
		\Delta=\sigma_{2}(P).
	\end{equation}
	For the purposes of large $n$ comparison, view $\|E\|_{2}\approx 2\sqrt{\Delta}$ from \cite{Lu-Peng--2013} and set the lower threshold $\alpha$ to be $\|E\|_{2}$. Define $r_{p,q}$ to be the edge probability-dependent parameter $r_{p,q}:=(p+q)/(p-q)$. Then via Kato--Temple methodology applied jointly to $\sigma_{1}(P)$ and $\sigma_{2}(P)$, with probability approximately 0.99 when $t_{KT} \ge 2.55$, for each singular value, respectively,
	\begin{align*}
		-3t_{KT}
		&\le \hat{\sigma}_{1}(A) - \sigma_{1}(P)
		\le 4r_{p,q} + t_{KT},\\
		-t_{KT} &\le \hat{\sigma}_{2}(A) - \sigma_{2}(P)
		\le (8+6t)r_{p,q}+t_{KT}.
	\end{align*}
	By the same approach, the bounds obtained in  \cite{O-Vu-Wang--2014} are given by
	\begin{align*}
		-t_{OVW} &\le \hat{\sigma}_{1}(A) - \sigma_{1}(P)
		\le 8\sqrt{2}r_{p,q} + \sqrt{2}t_{OVW},\\
		-t_{OVW} &\le \hat{\sigma}_{2}(A) - \sigma_{2}(P)
		\le 8+\sqrt{2}t_{OVW}.
	\end{align*}
	Direct application of the results in \cite{O-Vu-Wang--2014} yields probability approximately at least 0.99 for $t_{OVW} \ge 11.6$, though it appears upon further inspection that this can be improved to, for example, $t_{OVW} \ge 5.6$. The above joint analysis demonstrates that our bounds are favorable for the pair $\{\hat{\sigma}_{1}(A),\sigma_{1}(P)\}$ whereas the opposite is true for the pair $\{\hat{\sigma}_{2}(A),\sigma_{2}(P)\}$.
	
	We emphasize that here the upper bounds are of primary importance and interest. Indeed, the $(C,c,\gamma)$ property allows for straightforward lower bounds to be obtained by epsilon net techniques together with the Courant--Fisher--Weyl min-max principle. For example, note that a single application of $(C,c,\gamma)$-concentration yields that
	$\hat{\sigma}_{2}(A) - \sigma_{2}(P) \ge -t$
	with probability at least $1-C\exp(-ct^{\gamma})$.
	
	Among the advantages of Kato--Temple methodology is the ability to, in certain cases, refine one's initial analysis by further localizing the underlying interval $(\alpha,\beta)$. This is possible in the current example wherein we can ``zoom in'' further on the largest signal singular value. In particular, keeping the same indexing as above and setting $\alpha$ to be $\|E\|_{2}+\sigma_{1}(P)$, then for $n$ large and with probability approximately 0.99, we have
	\begin{flalign*}
		&&-t_{KT} &\le \hat{\sigma}_{2}(A) - \sigma_{2}(P)
		\le 2\left(\tfrac{p}{q}+1\right)+t_{KT}.
		&& \blacktriangle
	\end{flalign*}
	%Throughout this example, we note the bounds' dependence upon the underlying parameter $p$ and $q$. By virtue of the large $n$ comparison here, these parameters do not meaningfully influence the underlying probabilistic statement. \hfill$\blacktriangle$
\end{example}

In contrast to the low rank setting of Example \ref{ex:2BSSBM}, Example \ref{ex:SpikeModel} below demonstrates how our results can be applied to the problem of estimating signal in a high rank matrix setting.

\begin{example}[Estimating signal in a high rank spike model]
	\label{ex:SpikeModel}
	Let	$m,n,p \in \N$ and set $q:=m+n+p$. Let $M\in\R^{q \times q}$ be full rank with singular values given by the set
	\begin{align*}
	\{
	\underbrace{1,\dots,1}_{m \textnormal{ times }},
	\underbrace{\kappa+1,\dots,\kappa+1}_{n \textnormal{ times }},
	\underbrace{\tau+\kappa+1,\dots,\tau+\kappa+1}_{p \textnormal{ times }}
	\},
	\end{align*}
	where $\tau,\kappa > 0$. By slight abuse of notation, denote the singular values of $M$ up to multiplicity by $\sigma_{1}:=1$, $\sigma_{2}:=\kappa+1$, and $\sigma_{3}:=\tau+\kappa+1$.
	
	Further suppose that $E\in\R^{q \times q}$ has entries which are independent, identically distributed standard normal random variables. It follows by Gaussian concentration that $E$ is $(C,c,\gamma)$-concentrated with parameters $C=2$, $c=\frac{1}{2}$, and $\gamma=2$, and so by an application of Lemma \ref{lem:CcgammaSpectralProb}
	for $\epsilon=4$, then
	\begin{align*}
		\mathbb{P}\left[\|E\|_{2}>6\sqrt{q}\right]
		&\le 2\exp\left(-\tfrac{1}{10}q\right).
	\end{align*}
	Define $\hat{M}:=M+E$ and organize the singular values of $\hat{M}$ in correspondence with the repeated singular values of $M$, namely write
	\begin{equation*}
	\left\{
	\{\hat{\sigma}_{1,i_{1}}\}_{i_{1}=1}^{m},
	\{\hat{\sigma}_{2,i_{2}}\}_{i_{2}=1}^{n},
	\{\hat{\sigma}_{3,i_{3}}\}_{i_{3}=1}^{p}
	\right\}.
	\end{equation*}
	Suppose that $\tau,\kappa > (2\times(6\sqrt{q})+1)$. Then we can use Weyl's inequality as a preliminary tool for selecting the threshold values $\alpha$ and $\beta$. In particular, such analysis yields that with high probability,
	\begin{align*}
		|\hat{\sigma}_{1,i_{1}}-1|\le\|E\|_{2}
		&\Longrightarrow
		& 0
		&\le \hat{\sigma}_{1,i_{1}}
		\le 6\sqrt{q}+1,\\
		|\hat{\sigma}_{2,i_{2}}-(\kappa+1)|\le\|E\|_{2}
		&\Longrightarrow
		& 6\sqrt{q}+2
		&< \hat{\sigma}_{2,i_{2}}
		< \tau+\kappa-6\sqrt{q},\\
		|\hat{\sigma}_{3,i_{3}}-(\tau+\kappa+1)|\le\|E\|_{2}
		&\Longrightarrow
		&\tau+\kappa+1-6\sqrt{q}
		&\le \hat{\sigma}_{3,i_{3}}
		\le \tau+\kappa+1+6\sqrt{q}.
	\end{align*}
	For the choices $\alpha=6\sqrt{q}+2$ and $\beta=\tau+\kappa-6\sqrt{q}$, observe that $\{\hat{\sigma}_{2,i_{2}}\}_{i_{2}=1}^{n} \subset (\alpha,\beta)\subset \R_{>0}$ while simultaneously $\{1,\kappa+1,\tau+\kappa+1\}\bigcap(\alpha,\beta)=\{\kappa+1\}$. In this setting our perturbation theorems apply for $\kappa$ sufficiently large. Namely, choosing $\delta \in (0,1]$ and setting $t=\Theta(\log^{\delta}q)$ yields that for each $k\in[n]$ there exist positive constants $c^{\prime}$ and $c^{\prime\prime}$ such that with high probability,
	\begin{align*}
		|\hat{\sigma}_{2,k}-\sigma_{2}|
		&\le c^{\prime}t+c^{\prime\prime}.
	\end{align*}
	To reiterate, this bound improves upon the bound implied by a na\"{i}ve, terminal application of Weyl's inequality.
	Moreover, Example \ref{ex:SpikeModel} demonstrates how Weyl's inequality may be invoked for the preliminary purpose of establishing threshold values when the paired singular values (eigenvalues) correspond to the same index after ordering. \hfill$\blacktriangle$
\end{example}
%%%%%%%%%%%%%%%%%%%%%%%%%%%%%%%%%%%%%%%%%%%%%%%%%%%
%%%%%%%%%%%%%%%%%%%%%%%%%%%%%%%%%%%%%%%%%%%%%%%%%%%
%%%%%%%%%%%%%%%%%%%%%%%%%%%%%%%%%%%%%%%%%%%%%%%%%%%
\section{Applications to graph inference}
\label{sec:Application}
%%%%%%%%%%%%%%%%%%%%%%%%%%%%%%%%%%%%%%%%%%%%%%%%%%%
\subsection{Methods of graph inference}
The field of statistical inference and modeling for graphs represents a burgeoning area of research with implications for the social and natural sciences among other disciplines \cite{Kolaczyk2009,Goldenberg2010}. Within the current body of research, the pursuit of identifying and studying community structure within real--world networks continues to receive widespread attention \cite{Newman2006,NewmanGirvan2004,BickelChen2009,Fortunato2010,VerzelenCastro2015,CastroVerzelen2014}. Still another area of investigation involves anomaly detection for time series of graphs by considering graph statistics such as the total degree, number of triangles, and various scan statistics \cite{Wang2014, RukhinDiss}. Here we apply our results to two such detection tasks.
%%%%%%%%%%%%%%%%%%%%%%%%%%%%%%%%%%%%%%%%%%%%%%%%%%%
\subsection{Community detection via hypothesis testing}
\label{subsec:threeBlockSBM}
In this application we view the problem of community detection through the lens of hypothesis testing as in \cite{CastroVerzelen2014,VerzelenCastro2015}. We consider the simple setting of a balanced three block stochastic block model and the problem of detecting differences in between-block communication. Namely, consider the block edge probability matrix and block assignment vector given by
\begin{equation}
	\label{eqn:Null_B_0}
	\textnormal{Null model: }
	B_0 = \left( \begin{array}{ccc}
	p & q & q \\
	q & p & q \\
	q & q & p \end{array} \right)
	\text{ and }
	\pi_0 = \left(\tfrac{1}{3},\tfrac{1}{3},\tfrac{1}{3}\right),
\end{equation}
where $p = 0.81$ and $q = 0.2025$. In this model, vertices have an equal probability of belonging to each of the three blocks. Vertices within the same block have probability $p$ of being connected by an edge, whereas vertices in different blocks have probability $q$ of being connected by an edge.

As an aside, we note that this SBM may be cast in the language of random dot product graphs for which the underlying distribution of latent positions $F$ is a mixture of point masses. Specifically, take $F$ to be the discrete uniform distribution on the vectors $x_1 \approx (0.55, 0.32, 0.64)$, $x_2 \approx  (-0.55, 0.32, 0.64)$, and $x_3 \approx (0, -0.64, 0.64)$ in $\mathbb{R}^3$ (see Remarks \ref{rem:RDPG} and \ref{rem:kernel extension}).

For a graph on $n$ vertices from this three block model, condition on the graph exhibiting equal block sizes, i.e. $n_{1}=n_{2}=n_{3}=n/3$. For the corresponding $P$ matrix, denoted $P_n(B_0)$, the non-trivial (signal) model eigenvalues themselves exhibit multiplicity (hence Equation (\ref{eqn:CLT big 0_P bound}) via \cite{Athreya-et-al--2015} does not apply) and are
\begin{equation}
	\lambda_{1}(P_n(B_0))=\lambda_{2}(P_n(B_0))=\frac{n}{3}(p-q)
	\text{ and } \lambda_{3}(P_n(B_0))=\frac{n}{3}(p+2q).
\end{equation}

In contrast, consider an alternative model in which the first and second blocks exhibit stronger between-block communication. This stronger communication is represented by an additional additive factor $\epsilon\in(0,p-q)$ in the block edge probability matrix $B_{\epsilon}$, where $\epsilon$ is assumed to be bounded away from $p-q$ for convenience.
\begin{equation}
	\label{eqn:Alt_B_0}
	\textnormal{Alternative model: }
	B_{\epsilon} = \left( \begin{array}{ccc}
	p & q+\epsilon & q \\
	q+\epsilon & p & q \\
	q & q & p \end{array} \right)
	\text{ and }
	\pi_1 = \left(\tfrac{1}{3},\tfrac{1}{3},\tfrac{1}{3}\right).
\end{equation}
Under $B_{\epsilon}$, the signal eigenvalues of $P_{n}(B_{\epsilon})$ (equiv., singular values) can be explicitly computed as functions of $p,q,n,$ and $\epsilon$. They are given by \begin{align*}
	\lambda_{1}(P_{n}(B_{\epsilon}))
	&=\frac{n}{3}(p-q-\epsilon), \hspace{1em} \lambda_{2}(P_{n}(B_{\epsilon}))
	=\frac{n}{6}(2p+q+\epsilon-\sqrt{9q^2+2q \epsilon + \epsilon^2}),\\
	\lambda_{3}(P_{n}(B_{\epsilon}))
	&=\frac{n}{6}(2p+q+\epsilon+\sqrt{9q^2+2q \epsilon + \epsilon^2}).
\end{align*}
Furthermore, the maximum expected degree of the model corresponding to $B_{\epsilon}$ is given by $\Delta_{\epsilon}= \frac{n}{3}(p+2q+\epsilon)$.

For $\epsilon>0$, now consider a simple null versus simple alternative hypothesis test written as
\begin{equation}
	\label{eqn:hypothesis_test}
	\mathbb{H}_{0}: B = B_{0}
	\text{   against   }
	\mathbb{H}_{A}: B=B_{\epsilon}.
\end{equation}
In what follows we choose the smallest signal eigenvalue as our test statistic and denote it by $\Lambda_{1}$. We compare our bounds obtained via Kato--Temple methodology with the large-sample approximation bounds implied by \cite{Lu-Peng--2013} for the specified values $n\in\{6000,9000,12000,15000\}$. Similar comparison can be carried out with respect to the results in \cite{O-Vu-Wang--2014}. Our bounds compare favorably with those in \cite{O-Vu-Wang--2014} even for conservative choices of $t$ therein.

By Lemma \ref{lem:CcgammaSpectralProb} and Proposition \ref{prop:A-P}, irrespective of $\epsilon>0$ above, we have the concentration inequality $\mathbb{P}\left[\|E\|_{2}>3\sqrt{n}\right]
\le 2\exp\left(-\tfrac{1}{20}n\right)$.
This spectral norm bound allows us to invoke an unconditional version of Theorem \ref{thrm:IERMconditional}. Specifically, for moderate choices of $t>0$, the bounds in Theorem \ref{thrm:IERMconditional} hold with probability at least
$1-12\exp(-t^{2})-2\exp\left(-\frac{1}{20}n\right)$. When $n\ge 6000$, the choice $t \approx 2.66$ yields probability at least 0.99. 

Using these concentration inequality results, we determine confidence intervals which hold for $\Lambda_{1}$ with probability at least 0.99 under $\mathbb{H}_{0}$ and $\mathbb{H}_{A}$, respectively. We compute the value $\epsilon_{n}$ such that the confidence intervals under $\mathbb{H}_{0}$ and $\mathbb{H}_{A}$ no longer overlap for $\epsilon\in(\epsilon_{n},0.2]$, emphasizing that smaller values of $\epsilon_{n}$ indicate superior performance. This provides us with a region of the alternative in which the statistical test has power at least 0.99. Our results are summarized in the numerical table below.
\begin{table}[H]
	\caption{}
	%\begin{center}
	\begin{tabular}{||c c c||}
		\hline
		$n$ & $\epsilon_{n}$ via \cite{Lu-Peng--2013} & $\epsilon_{n}$ via this paper\\ [0.5ex] 
		\hline\hline
		6000 & 0.1006 & 0.0407 \\ 
		\hline
		9000 & 0.0818 & 0.0256 \\
		\hline
		12000 & 0.0707 & 0.0187 \\
		\hline
		15000 & 0.0631 & 0.0147 \\[1ex]
		\hline
	\end{tabular}
	%\end{center}
\end{table}
It is not too difficult to realize that the eigenvalue-based test considered here has asymptotic power equal to one as $n\rightarrow\infty$ for any choice of $0<q<p<1$ and $\epsilon\in(0,p-q)$. Moreover, as a consequence of Theorem \ref{thrm:spikeEvalsIERM} and subsequent discussion, we make the following observation.
\begin{proposition}
	Consider testing the hypothesis in Equation (\ref{eqn:hypothesis_test}). Assume that
	$q \equiv q_{n} = \omega(\frac{\log n}{n})$ with
	$q_{n} < p_{n}$. Then for $n \epsilon_{n} = \omega(\log n)$ and $\epsilon_{n} < p_{n} - q_{n}$, the above test using $\Lambda_{1}$ has asymptotically full power.
\end{proposition}

Note that the above analysis investigates testing performance as a function of $\epsilon$ for graphs with fixed block proportions. Next we investigate a setting wherein $\epsilon$ is fixed and the sizes of the graph communities change.

%%%%%%%%%%%%%%%%%%%%%%%%%%%%%%%%%%%%%%%%%%%%%%%%%%%
\subsection{Change-point detection}
We now consider a stylized example of change-point detection via hypothesis testing. Let $T^{*} \geq 1$ and suppose that $G_1, G_2, \dots, G_T$ for $T < T^{*}$ are Erd\H{o}s--R\'{e}nyi graphs on $n$ vertices, while for $T \ge T^{*}$ the graph $G_{T}$ is sampled according to a two block stochastic block model with block edge probability matrix $B = \Bigl[\begin{smallmatrix} p_{\epsilon} & p \\ p & p \end{smallmatrix} \Bigr]$ for $p_{\epsilon}:=p+\epsilon$ and $\epsilon>0$, with $m$ vertices assigned to the first block and $n - m$ vertices assigned to the second block. We note that $B$ encapsulates a notion of chatter anomaly, i.e., a subset of the vertices in $[n]$ exhibit altered communication behavior in an otherwise stationary setting. 
For a given value of $T$, we are interested in testing the hypothesis that $T$ is a change-point in the collection $\{G_1, G_2, \dots, G_T\}$. Given two graphs with adjacency matrices, $A^{(T-1)}$ and $A^{(T)}$, this can be formulated as the problem of testing the two-sample hypotheses
\begin{align*}
	\label{eq:two-sample}
	\mathbb{H}_0 &\colon A^{(T-1)} \sim \mathrm{ER}(n,p), A^{(T)} \sim \mathrm{ER}(n,p) \quad \text{against}\\
	\mathbb{H}_A &\colon A^{(T-1)} \sim \mathrm{ER}(n,p), A^{(T)} \sim \mathrm{SBM}(B, m, n-m).
\end{align*} 
We emphasize that in the above formulation, the parameter $p$ in $\mathrm{ER}(n,p)$, the size $m$ of the chatter community, and the associated communication probability $p_{\epsilon}$ are generally assumed to be unknown. 

Many test statistics are available for this change-point detection problem, including those based on graph invariant statistics (such as number of edges or number of triangles) or those based on locality statistics (such as max degree or scan statistics). For a given graph with adjacency matrix $A$, let $N(i) = \{j \colon A_{i,j} = 1\}$ denote the collection of vertices adjacent to vertex $i$. Furthermore,
\begin{itemize}
	\item let $\mathcal{T}_k$ count the number of $k$-cliques in $A$ for $k \geq 2$;
	\item let $\delta(A) := \max_{i} \sum_j A_{i,j}$ be the max degree statistic of $A$;
	\item let $\Psi(A) := \max_i \sum_{j,k \in N(i)} A_{j,k}$ be the scan statistic of $A$.
\end{itemize}
We note that these test statistics are widely used in anomaly detection for time series of graphs; see \cite{priebe_enron,Wang2014,CastroVerzelen2014,ranshous} and the references therein for a survey of results and applications.

One can then show \cite{rukhin2,tang_tsp2013} that the test statistics based on $\mathcal{T}_2$ and $\mathcal{T}_3$ are consistent for the above hypothesis test when $m = \Omega(\sqrt{n})$. More precisely, under the null hypothesis, one has
\begin{equation*}
	\frac{\mathcal{T}_2(A^{(T)}) - \mathcal{T}_2(A^{(T-1)})}{n\sqrt{p(1 - p)}} \overset{d}{\longrightarrow} N(0,1); \quad
	\frac{\mathcal{T}_3(A^{(T)}) - \mathcal{T}_3(A^{(T-1)})}{n^2 p^2 \sqrt{pp_{\epsilon}}} \overset{d}{\longrightarrow} N(0,1),
\end{equation*}
as $n \rightarrow \infty$,
while under the alternative hypothesis, one has
\begin{align*}
	\frac{\mathcal{T}_2(A^{(T)}) - \mathcal{T}_2(A^{(T-1)})}{n\sqrt{p(1 - p)}} &\overset{d}{\longrightarrow} N\Bigl(\frac{m(m-1)\epsilon}{n\sqrt{p(1-p)}}, C_1\Bigr); \\
	\frac{\mathcal{T}_3(A^{(T)}) - \mathcal{T}_3(A^{(T-1)})}{n^2 p^2 \sqrt{p p_{\epsilon}}}
	&\overset{d}{\longrightarrow} N\Bigl(\frac{\mu_{n,m,p,\epsilon}}{n^2p^2 \sqrt{pp_{\epsilon}}}, C_2\Bigr),
\end{align*}
as $n \rightarrow \infty$ for some positive constants $C_1$ and $C_2$ together with $\mu_{n,m,p,\epsilon}:=m^3p_{\epsilon}^{3}/6 + m^2(n-m)p^2p_{\epsilon} + (m(n-m)^2/2 + (n-m)^3/6)p^3 - n^3p^3/6$.
Now, if $m = \omega(\sqrt{n})$, then 
\begin{equation*}
	\frac{m(m-1)\epsilon}{n\sqrt{p(1-p)}} \rightarrow \infty;
	\quad \frac{\mu_{n,m,p,\epsilon}}{n^2p^2 \sqrt{pp_{\epsilon}}} \rightarrow \infty,
\end{equation*}
as $n \rightarrow \infty$,
and thus both $\mathcal{T}_2$ and $\mathcal{T}_3$ are consistent for the above hypothesis test when $m = \Omega(\sqrt{n})$. Furthermore, Theorem~2 and Proposition~2 of \cite{CastroVerzelen2014} indicate that $\mathcal{T}_2$ is asymptotically optimal, i.e., if $m = o(\sqrt{n})$ then provided that 
\begin{equation}
	\lim_{n \rightarrow \infty} \mathcal{I}(m,n,p,\epsilon)
	:= \lim_{n \rightarrow \infty} \frac{m \bigl(p_{\epsilon} \log \tfrac{p_{\epsilon}}{p} + (1 - p_{\epsilon}) \log \tfrac{1 - p_{\epsilon}}{1 - p}\bigr)}{2 \log{(n/m)}} < 1,
\end{equation}
no test statistic is consistent for testing the above hypotheses.
Similarly, one can also show \cite{rukhin1,tang_tsp2013} that the test statistics based on $\delta(A)$ and $\Psi(A)$ are consistent for the above hypothesis test when $m = \Omega(\sqrt{n \log n})$; in particular, the (normalized) limiting distributions of both $\delta(A^{(T)}) - \delta(A^{(T-1)})$ and $\Psi(A^{(T)}) - \Psi(A^{(T-1)})$ is the Gumbel distribution. 

In the context of this paper, one could also use a test statistic based on the largest eigenvalue. Our earlier results indicate that, under the null hypothesis, with high probability the largest eigenvalues of $A^{(T)}$ and $A^{(T-1)}$ satisfy
\begin{equation*}
	|\lambda_{\text{max}}(A^{(T)}) - \lambda_{\text{max}}(P^{(T)})| = O(1)
	\text{ and }|\lambda_{\text{max}}(A^{(T-1)}) - \lambda_{\text{max}}(P^{(T-1)})| = O(1),
\end{equation*}
along with $|\lambda_{\text{max}}(A^{(T)}) - \lambda_{\text{max}}(A^{(T-1)})| = O(1)$.
Meanwhile, under the alternative hypothesis, when $m=o(n)$, then with high probability
\begin{equation*}
	\begin{split}
		\left|\left|\lambda_{\text{max}}(A^{(T)}) - \lambda_{\text{max}}(A^{(T-1)})\right|
		- \frac{m^{2}p\epsilon}{np-m\epsilon}\right| &= O(1).
	\end{split}
\end{equation*}
Thus the largest eigenvalue test statistic is also consistent when $m = \Omega(\sqrt{n})$.

The previous test statistics are all global test statistics in the sense that, if $\mathbb{H}_0$ is rejected, the resulting test procedures do not extract the subset of the vertices which exhibits anomalous behavior between $A^{(T)}$ and $A^{(T-1)}$. One can construct related local test statistics which do extract the subset of anomalous vertices, although the resulting test procedure is computationally prohibitive. For example, assuming that $m$ is known, we could replace $\Psi(A)$ with the (modified) scan statistic $\Upsilon_m(A) = \max_{|S| = m} \mathcal{T}_{2}(A_{\mid S})$ where $A_{\mid S}$ is the subgraph of $A$ induced by the vertices in $S$ and the maximum is taken over all subsets $S \subset [n]$ with $|S| = m$. Thus $\Upsilon_m(A)$ is the maximum number of edges in any subgraph induced by $m$ vertices of $A$. By \cite{CastroVerzelen2014} the test statistic $\Upsilon_m(A^{(T)}) - \Upsilon_m(A^{(T-1)})$ is consistent for the hypothesis test considered in this section whenever 
\begin{equation*}
	\lim_{n \rightarrow \infty} \mathcal{I}(m,n,p,\epsilon) > 1.
\end{equation*}
Thus, for any fixed $p$ and $\epsilon$, the (modified) scan statistic is consistent when $m = \Omega(\log{n})$ as $n \rightarrow \infty$.
Using a similar idea, one can define a local variant of the largest eigenvalue statistic as $\Lambda_m(A) = \max_{|S| = m} \lambda_{\text{max}}(A_{\mid S})$. By Theorem~\ref{thrm:IERMconditional} and a union bound over all $\tbinom{n}{m} =  O(n^m)$ subsets $S \subseteq [n]$ with $|S| = m$, we have that there exists a constant $C > 0$ such that if $T = C \sqrt{m \log n}$, then with high probability
\begin{equation*}
	|\Lambda_m(A^{(T)}) - \Lambda_m(A^{(T-1)})| = O(\sqrt{m \log n})
\end{equation*}
under the null hypothesis, whereas under the alternative hypothesis, with high probability
\begin{equation*}
	\left| |\Lambda_m(A^{(T)}) - \Lambda_m(A^{(T-1)})| - m\epsilon\right|
	= O(\sqrt{m \log n}).
\end{equation*}
Thus for any fixed $p$ and $\epsilon$, the test statistic based on $\Lambda_m$ is also consistent for the above hypothesis test whenever $m = \Omega(\log{n})$ as $n \rightarrow \infty$.

In summary, the results in Section \ref{sec:Results} facilitate eigenvalue-based test statistics for the change-point detection problem as presented in this section. Furthermore, the resulting procedure is consistent whenever the size of the chatter community $m$ exceeds the threshold of detectability given in \cite{CastroVerzelen2014}.
%%%%%%%%%%%%%%%%%%%%%%%%%%%%%%%%%%%%%%%%%%%%%%%%%%%
%%%%%%%%%%%%%%%%%%%%%%%%%%%%%%%%%%%%%%%%%%%%%%%%%%%
%%%%%%%%%%%%%%%%%%%%%%%%%%%%%%%%%%%%%%%%%%%%%%%%%%%
\section{Acknowledgments}
\label{sec:Thanks}
The authors thank the anonymous referees for their
valuable feedback which improved the quality of this paper.
%%%%%%%%%%%%%%%%%%%%%%%%%%%%%%%%%%%%%%%%%%%%%%%%%%%
%%%%%%%%%%%%%%%%%%%%%%%%%%%%%%%%%%%%%%%%%%%%%%%%%%%
%%%%%%%%%%%%%%%%%%%%%%%%%%%%%%%%%%%%%%%%%%%%%%%%%%%
\section{Appendix}
\label{sec:Appendix}
%%%%%%%%%%%%%%%%%%%%%%%%%%%%%%%%%%%%%%%%%%%%%%%%%%%
\subsection{Proof of Theorem \ref{thrm:IERMconditional}}
\begin{proof}
	Let $P,E \in \R^{n \times n}$ be real symmetric matrices such that $E$ satisfies Proposition \ref{prop:A-P}.
	%is $(C,c,\gamma)$ concentrated and the perturbation of $P$ by $E$ is defined as $A:=P+E$.
	Denote the $d$ largest eigenvalues of $P$ and $A$ by
	\begin{align*}
		0 &< \lambda_{1}(P) \le \lambda_{2}(P) \le \dots \le \lambda_{d}(P),\\
		0 &< \lambda_{1}(A) \le \lambda_{2}(A) \le \dots \le \lambda_{d}(A).
	\end{align*}
	Let $\{w_{i}\}_{i=1}^{d}$ denote a collection of orthonormal eigenvectors of $P$ corresponding to the collection of eigenvalues $\{\lambda_{i}(P)\}_{i=1}^{d}$. Similarly, let $\{u_{i}\}_{i=1}^{d}$ denote a collection of orthonormal eigenvectors of $A$ corresponding to the collection of eigenvalues $\{\lambda_{i}(A)\}_{i=1}^{d}$.
	
	For each $i \in [d]$ define $\eta_{i}$ to be an ``approximate eigenvalue of $A$ close to $\lambda_{i}(P)$'' in the sense that
	\begin{align}
		\eta_{i} &:= \langle Aw_{i},w_{i}\rangle = \lambda_{i}(P) + \langle Ew_{i},w_{i}\rangle,
	\end{align}
	and define a corresponding ``residual quantity'' $\epsilon_{i}$ as
	\begin{align}
		\epsilon_{i}:=\|(A-\eta_{i})w_{i}\|.
	\end{align}
%%%%%%%%%%%%%%%%%%%%%%%%%%%%%%%%%%%%%%%%%%%%%%%%%%%
\subsubsection{Proof of Theorem \ref{thrm:IERMconditional}: upper bound}
	\label{sec:Upper bound}	
	Now for fixed $k\in[d]$ define the $k$-dimensional linear manifold $\mathcal{M}_{k}$ by
	\begin{align*}
		\mathcal{M}_{k}:=\textnormal{span}\{u_{1},\dots,u_{k}\}.
	\end{align*}
	We now define a collection of ``aggregate quantities'':
	\begin{itemize}
		\item Define $w$ to be an ``aggregate approximate eigenvector of $A$'' in the sense that
		$ w:=\sum_{i=1}^{k}r_{i}w_{i} $
		for a collection of normalized coefficients $\{r_{i}\}_{i=1}^{k}$ such that
		$ \|w\|^{2} = \sum_{i=1}^{k}r_{i}w_{i} = 1,$
		and satisfying the under-determined linear system
		$ \langle w, u_{i} \rangle = 0 \textnormal{ for } i=1,2,\dots,k-1.$
		\item Define $\eta$ to be an ``aggregate approximate eigenvector of $A$'' in the sense that $\eta := \langle Aw,w \rangle.$
		\item Define $\epsilon$ to be the ``aggregate residual quantity'' $\epsilon:=\|(A-\eta)w\|.$
	\end{itemize}
	By Lemma 1 in \cite{Kato--1950}, the interval $\left(\alpha, \eta + \frac{\epsilon^2}{\eta - \alpha}\right]$ contains a point in the spectrum of $A$. Note that by construction, $w \in \mathcal{M}_{k-1}^{\perp} = : \mathcal{N}_{k-1}$; moreover, $Aw \in \mathcal{N}_{k-1}$ as a function of $\{r_{i}\}_{i=1}^{k}$. In the Hilbert space $\mathcal{N}_{k-1}$, however, the spectrum of $A$ does not contain $\lambda_{1}(A), \dots, \lambda_{k-1}(A)$ since $u_{1},\dots,u_{k-1} \notin \mathcal{N}_{k-1}$. Thus, by another  application of Lemma 1 in \cite{Kato--1950}, the eigenvalue of $A$ in the interval given by $\left(\alpha, \eta + \frac{\epsilon^2}{\eta - \alpha}\right]$ must be $\lambda_{k}(A)$ with associated unit eigenvector $u_{k}$. Hence,
	\begin{equation}
		\label{eqn:Kato general upper bound}	
		\lambda_{k}(A)
		\le \eta + \frac{\epsilon^2}{\eta - \alpha}
		= \frac{\eta^2+\epsilon^2-\alpha \eta}{\eta - \alpha}.
	\end{equation}
	We pause briefly to make several computational observations.
	First,
	\begin{align}
		\eta^2 + \epsilon^2
		&= \langle Aw,w \rangle^2 + \|(A-\langle Aw,w\rangle)w\|^2\\
		&= \|Aw\|^2
		=\sum_{i,j=1}^{k}r_{i}r_{j}\langle Aw_{i},Aw_{j}\rangle.	
	\end{align}
	Letting $\delta_{i,j}:=\mathbb{I}\{i=j\}$ denote the Kronecker delta function, we have for each $i,j\in[d]$ that
	\begin{align}
		\label{eqn:Awi ip Awj}
		\langle Aw_{i},Aw_{j} \rangle
		= \langle (A-\eta_{i})w_{i}, (A-\eta_{j})w_{j} \rangle + (\eta_{i} + \eta_{j} ) \langle Aw_{i},w_{j}\rangle - \eta_{i}^2 \delta_{i,j},
	\end{align}	
	while
	\begin{align}
		\label{eqn:Awi,wj,iNotj}
		\langle Aw_{i},w_{j} \rangle
		&= \langle Ew_{i},w_{j} \rangle
		\textnormal{ for } i \neq j.
	\end{align}
	It will also prove useful to recognize the expansion
	\begin{align}
		\eta
		&= \langle Aw,w \rangle
		= \sum_{i=1}^{k}r_{i}^{2}\eta_{i} + \sum_{1\le i < j \le k}2r_{i}r_{j}\langle Ew_{i},w_{j}\rangle.
	\end{align}
	Combining these observations yields
	\begin{align*}
		\eta^2 + \epsilon^2
		&= \sum_{i,j=1}^{k}r_{i}r_{j} \langle Aw_{i},Aw_{j} \rangle \\
		&= \left( \sum_{i=1}^{k}r_i^2 \langle Aw_{i},Aw_{i} \rangle \right) + \left( \sum_{1\le i < j \le k}2r_{i}r_{j}\langle Aw_{i},Aw_{j} \rangle \right) \\
		&=  \sum_{i,j=1}^{k}r_{i}r_{j} \langle (A-\eta_{i})w_{i}, (A-\eta_{j})w_{j} \rangle \\
		&+ \sum_{i=1}^{k}r_{i}^2\eta_{i}^2
		+ \sum_{1\le i<j \le k}2r_{i}r_{j} (\eta_{i} + \eta_{j} ) \langle Aw_{i},w_{j} \rangle.
	\end{align*}
	An application of the Cauchy--Schwarz inequality coupled with subsequent computation yields
	\begin{align*}
		\sum_{i,j=1}^{k}r_{i}r_{j}\langle(A-\eta_{i})w_{i},(A-\eta_{j})w_{j} \rangle
		&\le \sum_{i,j=1}^{k}r_{i}r_{j}\left(\|(A-\eta_{i})w_{i}\| \|(A-\eta_{j})w_{j}\|\right)\\
		&= \sum_{i,j=1}^{k}(r_{i}\epsilon_{i})(r_{j}\epsilon_{j})
		\le \left(\sum_{i=1}^{k}\epsilon_{i}|r_{i}|\right)^2.
	\end{align*}	
	Hence,
	\begin{align}
		\label{eqn:eta2epsilon2}
		\eta^{2} + \epsilon^{2}
		&\le \left(\sum_{i=1}^{k}\epsilon_{i}|r_{i}|\right)^2
		+ \sum_{i=1}^{k}r_{i}^2\eta_{i}^2
		+ \sum_{1\le i<j \le k}2r_{i}r_{j} (\eta_{i} + \eta_{j} ) \langle Ew_{i},w_{j} \rangle.
	\end{align}
	Returning to Eqn. (\ref{eqn:Kato general upper bound}), the numerator then becomes
	\begin{equation}	
		\left(\sum_{i=1}^{k}\epsilon_{i}|r_{i}|\right)^2
			+ \sum_{i=1}^{k}r_{i}^{2}\eta_{i}(\eta_{i} - \alpha)
			+ \sum_{1\le i<j \le k}2r_{i}r_{j}(\eta_{i} + \eta_{j} - \alpha)\langle Ew_{i},w_{j} \rangle
	\end{equation}
	while the denominator becomes
	\begin{equation}
		\left(\sum_{i=1}^{k}r_{i}^{2}(\eta_i - \alpha)\right) +
		\left( \sum_{1\le i < j \le k}2r_{i}r_{j}\langle Ew_{i},w_{j} \rangle \right).
	\end{equation}	
	By a simple union bound, observe that for $t>0$,
	\begin{align}
		\mathbb{P}\left[\textnormal{max}_{1 \le i \le j \le k}|\langle Ew_{i},w_{j}\rangle|>t\right]
		&\le \left(k+\binom{k}{2}\right)C\exp(-ct^{\gamma}),
	\end{align}
	in which case with high probability,
	\begin{align}
		\left(\sum_{i=1}^{k}r_{i}^{2}(\eta_i - \alpha)\right)
		&\ge \lambda_{1}(P)-\alpha-t,\\
		\left( \sum_{1\le i < j \le k}2r_{i}r_{j}\langle Ew_{i},w_{j} \rangle \right)
		&\ge -k(k-1)t,
	\end{align}
	while with high probability,
	\begin{align}
		\sum_{1\le i<j \le k}2r_{i}r_{j}(\eta_{i} + \eta_{j} - \alpha)\langle Ew_{i},w_{j} \rangle
		&\le (2\lambda_{k}(P)-\alpha +2t)k(k-1)t,\\
		\left(\sum_{i=1}^{k}r_{i}^{2}\eta_{i}\right)k(k-1)t
		&\le (\lambda_{k}(P)+t)k(k-1)t.
	\end{align}
	By adding and subtracting $\left(\sum_{i=1}^{k}r_{i}^{2}\eta_{i}\right)k(k-1)t$ to the numerator of Eqn. (\ref{eqn:Kato general upper bound}) we obtain the following bound in which the first term on the right-hand side is the leading term while the second term on the right hand side corresponds to a residual term.
	\begin{align*}
		\lambda_{k}(A)
		&\le \frac{ \left( \sum_{i=1}^{k}\epsilon_i |r_i| \right)^2 + \left(\sum_{i=1}^{k}r_i^2 \eta_i (\eta_i - \alpha -k(k-1) t) \right) }{\left(\sum_{i=1}^{k}r_i^2(\eta_i - \alpha - k(k-1) t)  \right)} \\
		&+ \frac{(3\lambda_{k}(P)-\alpha +3t)k(k-1)t}{\lambda_{1}(P)-\alpha - (k(k-1)+1)t}.
	\end{align*}
	Now by the same arguments as in \cite{Kato--1950}, Section 3, Eqns. (22--30), the constants $\{r_{i}\}_{i=1}^{k}$ can be removed. To this end, the quantity
	\begin{align}
		\frac{ \left( \sum_{i=1}^{k}\epsilon_i |r_i| \right)^2 + \left(\sum_{i=1}^{k}r_i^2 \eta_i (\eta_i - \alpha -k(k-1) t) \right) }{\left(\sum_{i=1}^{k}r_i^2(\eta_i - \alpha - k(k-1) t)  \right)}
	\end{align}
	is bounded above by the quantity
	\begin{align}
		\underset{1\le i \le k}{\textnormal{max}}\eta_{i} 
		+ \left(\sum_{i=1}^{k}\frac{\epsilon_{i}^2}{\eta_{i} - \alpha - k(k-1)t}\right).
	\end{align}
	Note that $\underset{1\le i \le k}{\textnormal{max}}\eta_{i} \le \lambda_{k}(P)+t,$ with high probability, while a simple computation reveals that for each $i\in[k]$,
	\begin{align}
		\label{eqn:epsilonSquared}
		\epsilon_{i}^{2}
		&= \|Ew_{i}\|^{2}-|\langle Ew_{i},w_{i}\rangle|^{2}
		\le \|Ew_{i}\|^{2}
		\le \|E\|_{2}^{2}.
	\end{align}
	Putting all these observations together finally produces an upper bound on $\lambda_{k}(A)$ of the form
	\begin{align}
		\lambda_{k}(A)
		\le \lambda_{k}(P) + t
		&+ \zeta^{+},
	\end{align}
	where $\zeta^{+}:=\frac{k\|E\|_{2}^{2}+(3\lambda_{k}(P)-\alpha +3t)k(k-1)t}{\lambda_{1}(P)-\alpha - (k(k-1)+1)t}$.
%%%%%%%%%%%%%%%%%%%%%%%%%%%%%%%%%%%%%%%%%%%%%%%%%%%
\subsubsection{Proof of Theorem \ref{thrm:IERMconditional}: lower bound}
	Fix $k\in[d]$ and let $l:=d-k+1$. Define $\mathcal{M}_{l}$ to be the $l$-dimensional linear manifold given by
	\begin{align*}
		\mathcal{M}_{l}:=\textnormal{span}\{u_{k},\dots,u_{d}\}.
	\end{align*}
	We now define a collection of ``aggregate quantities'' similar to the formulation in Section \ref{sec:Upper bound}:
	\begin{itemize}
		\item Define $w$ to be an ``aggregate approximate eigenvector of $A$'' in the sense that
		$ w:=\sum_{i=k}^{d}r_{i}w_{i} $
		for a collection of normalized coefficients $\{r_{i}\}_{i=k}^{d}$ such that
		$ \|w\|^{2} = \sum_{i=k}^{d}r_{i}w_{i} = 1,$
		and satisfying the under-determined linear system
		$ \langle w, u_{i} \rangle = 0 \textnormal{ for } i=k+1,\dots,d.$
		\item Define $\eta$ to be an ``aggregate approximate eigenvector of $A$'' in the sense that $\eta := \langle Aw,w \rangle$.
		\item Define $\epsilon$ to be the ``aggregate residual quantity'' $\epsilon:=\|(A-\eta)w\|.$
	\end{itemize}	
	By Lemma 2 in \cite{Kato--1950}, the interval $\left[\eta-\frac{\epsilon^2}{\beta-\eta}, \beta \right)$ contains a point in the spectrum of $A$. Note that by construction, $w\in \mathcal{M}_{l-1}^{\perp} = : \mathcal{N}_{l-1}$; moreover, $Aw \in \mathcal{N}_{l-1}$ as a function of $\{r_{i}\}_{i=k}^{d}$. In the Hilbert space $\mathcal{N}_{l-1}$, however, the spectrum of $A$ does not contain $\lambda_{k+1}(A), \dots, \lambda_{d}(A)$ since $u_{k+1},\dots,u_d \notin \mathcal{N}_{l-1}$. Thus, by another  application of Lemma 2 in \cite{Kato--1950}, the eigenvalue of $A$ in the interval $\left[\eta-\frac{\epsilon^2}{\beta-\eta}, \beta \right)$ must be $\lambda_{k}(A)$ with associated unit eigenvector $u_{k}$.
	
	Consider first the special case when $\beta=\infty$. By a simple union bound, observe that for $t>0$,
	\begin{align}
		\mathbb{P}\left[\textnormal{max}_{k \le i \le j \le d}|\langle Ew_{i},w_{j}\rangle|>t\right]
		&\le \left(l+\binom{l}{2}\right)C\exp(-ct^{\gamma}),
	\end{align}
	hence with high probability
	\begin{equation}
		\label{eqn:partial denom bound}
		\sum_{i=k}^{d} r_{i}^2 \eta_{i}
		\ge \underset{k \le i \le d}{\text{min}} \eta_i
		\ge \lambda_{k}(P)-t
	\end{equation}
	and
	\begin{align}
		\label{eqn:applicKatoLowerBound}
		\lambda_{k}(A)
		\ge \eta
		&= \sum_{i=k}^{d} r_i^2 \eta_i + \sum_{\substack{k\le i<j \le d}}2r_{i}r_{j} \langle Ew_{i},w_{j} \rangle\\
		&\ge \lambda_{k}(P) - (l(l-1)+1)t.
	\end{align}
	Now suppose that $\beta < \infty.$ Then for the lower bound of the above interval, one has
	\begin{align*}
	\lambda_{k}(A)
	&\ge \eta - \frac{\epsilon^2}{\beta- \eta}
	= \frac{\beta \eta-\eta^2-\epsilon^2}{\beta-\eta}
	=\frac{-(\eta^2+\epsilon^2)+\beta \eta}{\beta-\eta}.
	\end{align*}
	Reversing the direction of the previous application of the Cauchy--Schwarz inequality in Eqn. (\ref{eqn:eta2epsilon2}) permits the numerator to be bounded below by
	\begin{align*}
		-(\sum_{i=k}^{d}\epsilon_{i}|r_{i}|)^2
		+\sum_{i=k}^{d}r_{i}^{2}\eta_{i}(\beta-\eta_{i})
		+ \sum_{k\le i < j \le d}2r_{i}r_{j}(\beta- \eta_{i} -\eta_{j})\langle Ew_{i},w_{j}\rangle,
	\end{align*}
	whereas the denominator has the expansion
	\begin{align*}
		\sum_{i=k}^{d}r_{i}^{2}(\beta-\eta_{i})
		+ \sum_{k \le i < j \le d}2r_{i}r_{j}\langle Ew_{i},w_{j}\rangle.
	\end{align*}
	For the denominator terms, note that with high probability
	\begin{align*}
		\sum_{i=k}^{d}r_{i}^{2}(\beta-\eta_{i})
		&\ge \beta - \lambda_{d}(P) - t,\\
		\sum_{k \le i < j \le d}2r_{i}r_{j}\langle Ew_{i},w_{j}\rangle
		&\ge -l(l-1)t,
	\end{align*}
	while in the numerator, with high probability,
	\begin{align*}
		\sum_{k\le i < j \le d}2r_{i}r_{j}(\beta- \eta_{i} -\eta_{j})\langle Ew_{i},w_{j}\rangle
		&\ge -(\beta-\lambda_{k}(P)+\lambda_{d}(P)+2t)l(l-1)t.
	\end{align*}
	In the numerator of Eqn. (\ref{eqn:applicKatoLowerBound}), add and subtract the quantity $\left(\sum_{i=k}^{d}r_{i}^{2}\eta_{i}\right)l(l-1)t$
	which is bounded below by $(\lambda_{k}(P)-t)l(l-1)t$. Combining these observations yields
	%\hspace{1em}
	\begin{align*}
	\lambda_{k}(A)
	&\ge \frac{-(\sum_{i=k}^{d}\epsilon_i |r_i|)^2
		+\sum_{i=k}^{d}r_{i}^{2}\eta_{i}(\beta-\eta_{i} - l(l-1)t)}{\sum_{i=k}^{d}r_{i}^2(\beta-\eta_{i} - l(l-1)t)}\\
	&+ \frac{-(\beta-\lambda_{k}(P)+\lambda_{d}(P)-\lambda_{k}+3t)l(l-1)t}{\beta - \lambda_{d}(P) - (l(l-1)+1)t}.
	\end{align*}
	By employing the same approach used to obtain the upper bound and taking negatives when necessary (thereby reversing the direction in which bounds hold), we obtain the lower bound for $\lambda_{k}(A)$ of the form
	\begin{align}
		\lambda_{k}(A)
		\ge \lambda_{k}(P) - t
		&- \zeta^{-},
	\end{align}
	where $\zeta^{-} := \frac{l\|E\|_{2}^{2}+((\beta-\lambda_{k}(P))+(\lambda_{d}(P)-\lambda_{k}(P)) + 3t)l(l-1)t}{\beta-\lambda_{d}(P)-(l(l-1)+1)t}$.
\end{proof}%%%%%%%%%%%%%%%%%%%%%%%%%%%%%%%%%%%%%%%%

\subsection{Proof of Theorem \ref{thrm:spikeEvalsIERM}}
\begin{proof}
	The hypotheses imply by \cite{Lu-Peng--2013} that $\|A-P\|_{2}=O(\sqrt{\Delta})$ with probability $1-o(1)$ as $n\rightarrow\infty$. Set $\alpha=(C-c)\Delta/2$ and $\beta=\infty$ as Kato--Temple threshold values. Choose $\delta\in(0,1]$ and set $t=\Theta(\log^{\delta}n)$. Then in Theorem \ref{thrm:IERMconditional}, for sufficiently large $n$, one has $\zeta^{+},\zeta^{-}=O(t)$ where the underlying constant depends upon $k$, $d$, as well as underlying (unspecified) constants. So for $n\ge n_{0}$, then 
	$|\hat{\sigma}_{k}-\sigma_{k}| \le c_{k,d}t$ with probability $1-o(1)$ as claimed.
\end{proof}

\subsection{Proof of Theorem \ref{thrm:generalSVbound}}
\begin{proof}
	\label{pf:generalSVbound}
	The proof follows essentially \emph{mutatis mutandis} as in Theorem \ref{thrm:IERMconditional} via Remark \ref{rem:HermDial}, Definition \ref{def:Ccgamma}, and Lemma \ref{lem:CcgammaBlowUp}. In particular, observe that one has $\langle \tilde{\hat{M}}\tilde{w}_{i},\tilde{w}_{j}\rangle
	=\sigma_{i}\delta_{i,j}
	+\langle\tilde{E}\tilde{w}_{i},\tilde{w}_{j}\rangle$
	for each pair $i,j$, while at the same time
	$\|\tilde{E}\|_{2}=\|E\|_{2}$.
\end{proof}

\subsection{Proof of Lemma \ref{lem:CcgammaSpectralProb}}
\begin{proof}
Let $E\in\R^{m \times n}$ be a $(C,c,\gamma)$-concentrated random matrix. Take $\mathcal{X}$ and $\mathcal{Y}$ to be $\frac{1}{4}$-nets of the spheres $S^{n-1}$ and $S^{m-1}$, respectively, with cardinalities at most $9^{n}$ and $9^{m}$, respectively. Then a standard net argument yields that for $t>0$,
\begin{align*}
	\mathbb{P}\left[\|E\|_{2}>t\right]
	&\le \mathbb{P}\left[2 \underset{x\in\mathcal{X},y\in\mathcal{Y}}{\textnormal{max}}|\langle Ex,y\rangle| > t\right]\\
	&\le 9^{m+n}\mathbb{P}\left[|\langle Ex,y \rangle|>t/2\right]\\
	&\le C\exp((m+n)\log(9)-c(t/2)^{\gamma})\\
	&\le C\exp(2\log(9)\textnormal{max}\{m,n\}-c(t/2)^{\gamma}).
\end{align*}
Choose $\epsilon > 0$ such that $2 + \epsilon > 2\left(2\log(9)/c\right)^{1/\gamma}$ and set $t= (2+\epsilon)\textnormal{max}\{m,n\}^{1/\gamma}$. Then for $c_{\epsilon,c,\gamma}:=\left(c(1+\epsilon/2)^{\gamma}-2\log(9)\right)>0$, we have
\begin{align*}
	\mathbb{P}\left[\|E\|_{2}>(2+\epsilon)\textnormal{max}\{m,n\}^{1/\gamma}\right]
	&\le C\exp(-c_{\epsilon,c,\gamma}\textnormal{max}\{m,n\}).
\end{align*}
If in addition $m=n$ and $E$ is assumed to be symmetric, then since $\|E\|_{2}\equiv\textnormal{sup}_{\|x\|_{2}=1}|\langle Ex,x\rangle|$, one need only consider the $\frac{1}{4}$-net $\mathcal{X}$ for the purposes of a union bound.
\end{proof}

\newpage
\bibliographystyle{amsplain}
\bibliography{EJS1328}

\end{document}